\documentclass[a4paper, 12pt]{extarticle}
\usepackage{mathtext}
\usepackage{graphicx}
\usepackage{amscd}
\usepackage[T2A]{fontenc}
\usepackage[english]{babel}
\usepackage{latexsym,amsfonts,amssymb,amsmath,longtable,amsthm}
\usepackage[pic]{xy}
\usepackage{color}
\usepackage{bbm,makeidx}
\usepackage[centerlast,small]{caption}
\usepackage{amsmath}
\usepackage[cp1251]{inputenc}
\newtheorem{ttt}{{\scshape Theorem}}

\begin{document}
\title{Verbal quandles with one parameter}
\author{Elizaveta Markhinina, Timur Nasybullov}
\maketitle
\begin{abstract}
We find all words $W(x,y,z)$ in the free group $F(x,y,z)$, such that for every group $G$ and an element $c\in G$ the algebraic system $(G,*_{W,c})$ with the binary operation $*_{W,c}$ given by $a*_{W,c}b=W(a,b,c)$ for $a,b\in G$ is a quandle. Such quandles are called verbal quandles with one parameter.

~
\end{abstract}
\section{Introduction and preliminaries}
A quandle is an algebraic system with a binary operation whose three axioms encode the Reidemeister moves on knot diagrams. More formally, a quandle $Q$ is an algebraic system with one binary algebraic operation $(x,y)\mapsto x*y$ which satisfies the following three axioms:
\begin{itemize}
	\item[(q1)] $x*x=x$ for all $x\in Q$,
	\item[(q2)] the map $S_x:y\mapsto y*x$ is a bijection of $Q$ for all $x\in Q$,
	\item[(q3)] $(x*{y})*z=(x*z)*({y*z})$ for all $x,y,z\in Q$.
\end{itemize}

 First examples of quandles date back to 1940s when Takasaki defined keys, objects which were later known as involutory quandles \cite{Tak}. The concept of quandles was explicitly presented in the independent works of Joyce \cite{Joy} and Matveev \cite{Mat}. To each oriented diagram $D_K$ of an oriented knot $K$ in $\mathbb{R}^3$ they associate the quandle $Q(K)$ which does not change if we apply the Reidemeister moves to the diagram $D_K$. Joyce and Matveev proved that two knot quandles $Q(K_1)$ and $Q(K_2)$ are isomorphic if and only if there is a (probably, reversing orientation) homeomorphism of the ambient space $\mathbb{R}^3$ which maps $K_1$ to $K_2$. The knot quandle is a very strong invariant for knots in $\mathbb{R}^3$, however, usually it is very difficult
 to determine if two knot quandles are isomorphic. In \cite{brimi} it is shown that the isomorphism problem for quandles is, from the perspective of Borel reducibility, fundamentally difficult (Borel complete). Sometimes homomorphisms from knot quandles to simpler quandles provide useful information that helps determine whether two knot quandles  are isomorphic. This potential utility leads to the necessity of constructing quandles which are convenient to work with.

Another area where quandles find applications is the theory of the set-theoretical Yang-Baxter equation. The Yang-Baxter equation first appeared in theoretical physics and statistical
mechanics in the works of Yang \cite{Yan} and Baxter \cite{Bax1,Bax2}. The notion of the set-theoretical Yang-Baxter equation, was introduced by V.~Drinfel'd in the context of quantum groups (see \cite{Dri}). Recall that a set-theoretic solution of the Yang-Baxter equation is a pair $(X, r)$,
where $X$ is a set and $r:X\times X\to X\times X$ is a bijective map such that 
$$(r\times id)(id \times r)(r\times id)=(id\times r)(r \times id)(id\times r).$$
If $Q$ is a quandle with the operation $*$, then the pair $(Q,r)$, where $r:Q\times Q\to Q\times Q$ is the map given by 
\begin{equation}\label{ybeonquandle}
	r(x,y)=(y*x,x)
\end{equation}
for $x,y\in Q$ is a set-theoretic solution of the Yang-Baxter equation. 

In the recent papers \cite{BarNas1,BarNas3} there was constructed a general method how a given solution of the set-theoretical Yang-Baxter equation on an arbitrary algebraic system $X$ can be used for constructing a representation of the virtual braid group $VB_n$ by automorphisms of the algebraic system  $X$. In the papers \cite{byaka2,byaka1} these this method was used for constructing new representations of vertual braid groups. In the papers \cite{BarNas2, BarNas3}  there was introduced a method how a given solution of the set-theoretical Yang-Baxter equation on an algebraic system $X$ can be used for constructing an invariant for virtual knots and links which is an algebraic system from the same category as $X$. If $X$ is a quandle, and the solution of the set-theoretical Yang-Baxter equation on $X$ is given by (\ref{ybeonquandle}), then using procedures described in \cite{BarNas1,BarNas2,BarNas3} it is possible to construct a representation $VB_n\to {\rm Aut}(X)$, and a quandle invariant for virtual links. However, in order to apply procedures described in \cite{BarNas1,BarNas2,BarNas3} to $X$ and $r$, the quandle $X$ must be cinvenient to work with, hence it is necessary to construct such quandles.

Quandles were also investigated from an algebraic point of view and relations to other algebraic structures such as Lie algebras \cite{CarCra2}, Frobenius algebras and Yang-Baxter equation \cite{CarCra2}, Hopf algebras \cite{AndGra}, quasigroups and Moufang loops \cite{Moh1}, representation theory \cite{Moh2} and ring theory \cite{BarPas}. In all of these areas it is always useful to have examples of quandles which are convenient to work with.

A lot of examples of quandles which are convenient to work with come from groups. The most common example of a quandle is the conjugation quandle ${\rm Conj}(G)$ of a group $G$, i.~e. the quandle $(G,*)$ with $x*y=y^{-1}xy$ for $x,y\in G$. If we define another operation $*$ on the group $G$, namely $x*y=yx^{-1}y$, then the set $G$ with this operation also forms a quandle. This quandle is called the core quandle of a group $G$ and is denoted by ${\rm Core}(G)$. In particular, if $G$ is an abelian group, then the quandle ${\rm Core}(G)$ is called the Takasaki quandle of the abelian group $G$ and is denoted by $T(G)$. Such quandles were studied by Takasaki in \cite{Tak}. Other examples of quandles arising from groups can be found, for example, in \cite{BarDeySin, BarNasSin, BarNas6}.

The so-called verbal quandlуs were introduced in the paper \cite{BarNas4}. Let $w = w(x,y)$ be an element of the free group $F(x,y)$ with two generators. If $G$ is a group, then $w$ defines a map $w:G\times G\to G$ which maps a pair $(a,b)$ to $w(a,b)$. For $a,b\in G$ denote by $a *_w b = w(a, b)$. We can think about~$*_w$ as about new binary operation on the group $G$, so, $(G,*_w)$ is an algebraic system.  If the algebraic system $(G,*_w)$ is a quandle, then we call this quandle a verbal quandle defined by the word~$w$. It is clear that ${\rm Conj}(G)$ and ${\rm Core}(G)$ are verbal quandles. The following theorem which classifies verbal quandles is proved in \cite{BarNas4}.
\begin{ttt}\label{verbrack} Let $w\in F_2$ be an element from the free group. If $Q = (G, *_{w})$ is a quandle for every group $G$, then one of the following conditions holds
	\begin{enumerate}
		\item $w(x,y)=yx^{-1}y$,
		\item $w(x,y)=y^{-s}xy^{s}$ for $s\in\mathbb{Z}$.
	\end{enumerate}
\end{ttt}

In the present paper we introduce the notion of a verbal quandle with $n$ parameters which generalizes the notion of a verbal quandle. Let $W = W(x,y, z_1, z_2,\dots, z_{n})$ be a word from the free group with $(n+2)$ free generators. For a group $G$ and fixed elements $c_1,c_2,\dots,c_n\in G$ denote by $*_{W,c_1,c_2,\dots,c_n}$ the binary operation on $G$ defined by the rule 
$$a*_{W,c_1,c_2,\dots,c_n}b=W(a,b,c_1,c_2,\dots,c_n).$$
If $(G,*_{W,c_1,c_2,\dots,c_n})$ is a quandle, then this quandles is called a verbal quandle with $n$ parameters. If $n=0$, then this quandle is just a verbal quandle. The main result of the present paper is the following theorem which classifies verbal quandles with one parameter.

\begin{ttt}\label{th22233}	Let $W=W(x,y,z)$ be an element from the free group $F(x,y,z)$. If $Q = (G, *_{W,c})$ is a quandle for every group $G$ and an element $c\in G$, then one of the following conditions holds
	\begin{enumerate}
		\item $W(x,y,z)=yx^{-1}y$,
		\item $W(x,y,z)=y^{-s}xy^{s}$ for $s\in\mathbb{Z}$,
		\item $W(x,y,z)=yz^{-s}y^{-1}xz^{s}$ for $s\in \mathbb{Z}$,
		\item $W(x,y,z)=yz^{-s}xy^{-1}z^s$ for $s\in\mathbb{Z}$,
		\item $W(x,y,z)=z^{-s}y^{-1}xz^sy$ for $s\in\mathbb{Z}$,
		\item $W(x,y,z)=z^{-s}xy^{-1}z^sy$ for $s\in\mathbb{Z}$.
	\end{enumerate}
\end{ttt}
It is clear that Theorem~\ref{verbrack} follows from Theorem~\ref{th22233}. In the consequent paper we are going to apply solutions of the set-theoretic Yang-Baxter equation constructed by verbal quandles with one parameter using formula (\ref{ybeonquandle}) in order to construct new group invariants for virtual knots and links using constructions from \cite{BarNas1,BarNas2,BarNas3}.

~\\
\noindent\textbf{Acknowledgment.} The work is supported by Mathematical Center in Akademgorodok under agreement 075-15-2022-281 with the Ministry of Science and Higher Education of the Russian Federation.

\section{Proof of the main theorem}
Let $W= x^{\alpha_1} y^{\beta_1} z^{\gamma_1}x^{\alpha_2} y^{\beta_2} z^{\gamma_2}\ldots x^{\alpha_k} y^{\beta_k}z^{\gamma_k}$ for $\alpha_i, \beta_i,\gamma_i \in \mathbb{Z}$ be a reduced word in the free group $F(x,y,z)$. Since $(G,*_{W,c})$ is a quandle, from the axiom (q2) of a quandle follows that for every $a, b \in G$ there exists an element $d \in G$ such that 
\begin{equation}\label{start}d *_{W,c} a =W(d,a,c)=d^{\alpha_1} a^{\beta_1} c^{\gamma_1}d^{\alpha_2} a^{\beta_2} c^{\gamma_2}\ldots d^{\alpha_k} a^{\beta_k}c^{\gamma_k}= b.
\end{equation}
It meant that we can express $d$ from equation (\ref{start}) via $a,b,c$. It is possible if and only if $W =  u(y,z)  x^{\varepsilon} w(y,z)$ for $u(y,z),w(y,z)\in F(y,z)$, $\varepsilon \in \{\pm 1 \}$. From the axiom (q1) it follows that for every element $a\in G$ the equality $a*_{W,c}a=a$ holds. This equality can be rewritten in details in the following form 
$$a=a*_{W,c}a=W(a,a,c)=u(a,c)a^{\varepsilon}w(a,c),$$
and therefore $u(a,c)=a(w(a,c))^{-1}a^{-\varepsilon}$ for every group $G$, fixed element $c\in G$ and element $a\in G$. Hence we have $u(y,z)=y(w(y,z))^{-1}y^{-\varepsilon}$ and
\begin{equation}\label{generalequation}
W(x,y,z)=y(w(y,z))^{-1}y^{-\varepsilon}x^{\varepsilon}w(y,z).
\end{equation}

Let us understand what are the possible options for the word $w(y,z)\in F(y,z)$. In order to do this express the word $w(y,z)$ as a reduced word
 $$w = y^{q_1} z^{s_1}y^{q_2}z^{s_1} \ldots y^{q_n} y^{s_n},$$ 
 where $q_i, s_i \in \mathbb{Z}$ are non-zero integers with possible exceptions for $q_1$ and $s_n$. Denote by $|w(y,z)|$ the syllabic length  of the word $w(y,z)$. Depending on the word $w(y,z)$ and the parameter $\varepsilon\in\{\pm1\}$ we consider the cases collected in the following table.
 \begin{center}
 \begin{tabular}[t]{|l|l|l|l|}
 	\hline
 $|w|$	& $w\in F(y,z)$ & $\varepsilon\in\{\pm1\}$& $q_i,s_i\in\mathbb{Z}$ \\	\hline
 	1: $|w|=0$ & $1$ & & \\\hline
 2: $|w|=1$  & 2.1: $y^q$& & \\\cline{2-4}
   	& 2.2: $z^s$& 2.2.1: $\varepsilon=1$& \\\cline{3-4}
   	& & 2.2.2: $\varepsilon=-1$& \\\hline
 	3: $|w|=2$ & 3.1: $y^qz^s$ & 3.1.1: $\varepsilon=1$ &3.1.1.1: $q=-1$ \\\cline{4-4}
 	& & &3.1.1.2: $q\neq-1$ \\\cline{3-4}
  & & 3.1.2: $\varepsilon=-1$ &3.1.2.1: $q=1$ \\\cline{4-4}
  & &  &3.1.2.2: $q\neq1$ \\\cline{2-4}
  & 3.2: $z^sy^q$ & 3.2.1: $\varepsilon=1$ &3.2.1.1: $q=1$ \\\cline{4-4}
  & &  &3.2.1.2: $q\neq1$ \\\cline{3-4}
   & & 3.2.2: $\varepsilon=-1$ & \\\hline
 	4: $|w|=3$ & 4.1: $y^{q_1}z^sy^{q_2}$ & 4.1.1: $\varepsilon=1$&4.1.1.1: $q_2=1$, $q_1=-1$ \\\cline{4-4}
 	&& &4.1.1.2: $q_2=1$, $q_1\neq-1$ \\\cline{4-4}
 	&& &4.1.1.3: $q_2\neq1$ \\\cline{3-4}
 && 4.1.2: $\varepsilon=-1$&4.1.2.1: $q_2=1$ \\\cline{4-4}
 && &4.1.2.2: $q_2\neq1$ \\\cline{2-4}
  & 4.2: $z^{s_1}y^qz^{s_2}$ & &4.2.1: $|q|\neq 1$ \\\cline{4-4}
  & & &4.2.2: $|q|=1$ \\\hline
 	5: $|w|\geq 4$ & & &\\\hline
 \end{tabular}
\end{center}
From this table it is clear that the cases we are going to consider are all possible cases (the empty fields mean that the corresponding parameters can be arbitrary). We are going to prove that for the word $W(x,y,z)$ given by formula (\ref{generalequation}) the algebraic system $Q = (G, *_{W,c})$ is a quandle for every group $G$ and an element $c\in G$ if and only if the word $w(y,z)$ in formula (\ref{generalequation}) is given by one of the cases 1, 2.1, 2.2.1, 3.1.1.1, 3.2.1.1, 4.1.1.1. These cases give six possibilities for the word $W(x,y,z)$ described in the formulation of the theorem.

 \textbf{Case 1:} $|w(y,z)|=0$. In this case $w(y,z)=1$, and the word $$W(x,y,z)=y^{1-\varepsilon}x^{\varepsilon}$$ 
 doesn't depend on $z$, hence for every group $G$ we have $(G,*_{W,c})=(G,*_{W})$. From Theorem~\ref{verbrack} it follows that we have either $W(x,y,z)=yx^{-1}y$, or $W(x,y,z)=y^{-s}xy^{s}$ for $s\in \mathbb{Z}$. It means that in Case 1 we have a quandle.
 
 \textbf{Case 2:} $|w(y,z)|=1$. In this case we have either $w(y,z)=y^q$ for non-zero $q\in\mathbb{Z}$, or $w(y,z)=z^s$ for non-zero $s\in\mathbb{Z}$.

 If $w(y,z)=y^q$ for non-zero $q\in\mathbb{Z}$ (case 2.1), then the word $$W(x,y,z)=y^{1-q-\varepsilon}x^{\varepsilon}y^q$$ 
 doesn't depend on $z$, hence for every group $G$ we have $(G,*_{W,c})=(G,*_{W})$. From Theorem~\ref{verbrack} it follows that we have either $W(x,y,z)=yx^{-1}y$, or $W(x,y,z)=y^{-s}xy^{s}$ for $s\in \mathbb{Z}$. It means that in Case 2.1 we have a quandle.
 
  If $w(y,z)=z^s$ for non-zero $s\in\mathbb{Z}$ (case 2.2), then $W(x,y,z)=yz^{-s}y^{-\varepsilon}x^{\varepsilon}z^s$. Let us understand when the operation $*_{W,c}$ satisfies the third quandle axiom (q3) for each group $G$ and an element $c\in G$
  \begin{align}
  	\notag(a*_{W,c}b)*_{W,c}d&=(W(a,b,c))*_{W,c}d\\
  	\notag&=(bc^{-s}b^{-\varepsilon}a^{\varepsilon}c^s)*_{W,c}d\\
  	\notag&=W(bc^{-s}b^{-\varepsilon}a^{\varepsilon}c^s,d,c)\\
  	\label{left1}&=dc^{-s}d^{-\varepsilon}(bc^{-s}b^{-\varepsilon}a^{\varepsilon}c^s)^{\varepsilon}c^s,\\
 \notag 	(a*_{W,c}d)*_{W,c}(b*_{W,c}d)&=(W(a,d,c))*_{W,c}(W(b,d,c))\\
\notag&=(dc^{-s}d^{-\varepsilon}a^{\varepsilon}c^s)*_{W,c}(dc^{-s}d^{-\varepsilon}b^{\varepsilon}c^s)\\
\notag&=W(dc^{-s}d^{-\varepsilon}a^{\varepsilon}c^s,dc^{-s}d^{-\varepsilon}b^{\varepsilon}c^s,c)\\
\label{right1}&=(dc^{-s}d^{-\varepsilon}b^{\varepsilon}c^s)c^{-s}(dc^{-s}d^{-\varepsilon}b^{\varepsilon}c^s)^{-\varepsilon}(dc^{-s}d^{-\varepsilon}a^{\varepsilon}c^s)^{\varepsilon}c^s.
  \end{align}
If $\varepsilon=1$ (case 2.2.1), then from formulas (\ref{left1}), (\ref{right1}) it follows that 
\begin{align*}
(a*_{W,c}b)*_{W,c}d&=dc^{-s}d^{-1}bc^{-s}b^{-1}ac^{2s},\\
(a*_{W,c}d)*_{W,c}(b*_{W,c}d)&=dc^{-s}d^{-1}bc^sc^{-s}c^{-s}b^{-1}dc^sd^{-1}dc^{-s}d^{-1}ac^sc^s\\
&=dc^{-s}d^{-1}bc^{-s}b^{-1}ac^{2s},
\end{align*}
i.~e. the equality $(a*_{W,c}b)*_{W,c}d=(a*_{W,c}d)*_{W,c}(b*_{W,c}d)$ holds independently on $s$. It means that in Case 2.2.1 we have a quandle. If $\varepsilon=-1$ (case 2.2.2), then from formulas (\ref{left1}), (\ref{right1}) it follows that 
\begin{align}
\label{left2}	(a*_{W,c}b)*_{W,c}d&=dc^{-s}dc^{-s}ab^{-1}c^sb^{-1}c^s,\\
\notag	(a*_{W,c}d)*_{W,c}(b*_{W,c}d)&=dc^{-s}db^{-1}c^sc^{-s}dc^{-s}db^{-1}c^sc^{-s}ad^{-1}c^sd^{-1}c^s\\
\label{right2}	&=dc^{-s}db^{-1}dc^{-s}db^{-1}ad^{-1}c^sd^{-1}c^s.
\end{align}
Let us look to $a,b,c,d$ as to the generators of the free group $F(a,b,c,d)$ on four generators. In this situation from equalities (\ref{left2}), (\ref{right2}) it is clear that 
$$(a*_{W,c}b)*_{W,c}d\neq (a*_{W,c}d)*_{W,c}(b*_{W,c}d)$$
independently on $s$, i.~e. there exists a group $G=F(a,b,c,d)$ and an element $c\in G$ such that in $(G,*_{W,c})$ the third quandle axiom (q3) doen't hold, i.~e. $(G,*_{W,c})$ is not a quandle. It means that in Case 2.2.2 we don't have a quandle.

\textbf{Case 3:} $|w(y,z)|=2$. In this case $w(y,z)$ is either equal to $w(y,z)=y^qz^s$ or to $w(y,z)=z^sy^q$ for non-zero integers $q,s\in\mathbb{Z}$.

If $w(y,z)=y^qz^s$ (case 3.1), then $W(x,y,z)=yz^{-s}y^{-q-\varepsilon}x^{\varepsilon}y^qz^s$. Let us understand when the operation $*_{W,c}$ satisfies the third quandle axiom (q3) for each group $G$ and an element $c\in G$
  \begin{align}
	\notag(a*_{W,c}b)*_{W,c}d&=(W(a,b,c))*_{W,c}d\\
	\notag&=(bc^{-s}b^{-q-\varepsilon}a^{\varepsilon}b^qc^s)*_{W,c}d\\
	\notag&=W(bc^{-s}b^{-q-\varepsilon}a^{\varepsilon}b^qc^s,d,c)\\
	\label{left3}&=dc^{-s}d^{-q-\varepsilon}(bc^{-s}b^{-q-\varepsilon}a^{\varepsilon}b^qc^s)^{\varepsilon}d^qc^s,\\
	\notag 	(a*_{W,c}d)*_{W,c}(b*_{W,c}d)&=(W(a,d,c))*_{W,c}(W(b,d,c))\\
	\notag&=(dc^{-s}d^{-q-\varepsilon}a^{\varepsilon}d^qc^s)*_{W,c}(dc^{-s}d^{-q-\varepsilon}b^{\varepsilon}d^qc^s)\\
	\notag&=W(dc^{-s}d^{-q-\varepsilon}a^{\varepsilon}d^qc^s,dc^{-s}d^{-q-\varepsilon}b^{\varepsilon}d^qc^s,c)\\
	\notag&=(dc^{-s}d^{-q-\varepsilon}b^{\varepsilon}d^qc^s)c^{-s}(dc^{-s}d^{-q-\varepsilon}b^{\varepsilon}d^qc^s)^{-q-\varepsilon}\\
	\label{right3}&~~~~\cdot(dc^{-s}d^{-q-\varepsilon}a^{\varepsilon}d^qc^s)^{\varepsilon}(dc^{-s}d^{-q-\varepsilon}b^{\varepsilon}d^qc^s)^qc^s.
\end{align}
Let us consider the case when $\varepsilon=1$ (case 3.1.1). From equalities (\ref{left3}), (\ref{right3}) we see that the equality $(a*_{W,c}b)*_{W,c}d=(a*_{W,c}d)*_{W,c}(b*_{W,c}d)$ holds if an only if 
\begin{equation}\label{equality2}
c^{-s}b^{-q-1}ab^qc^sd^q=d^q(dc^{-s}d^{-q-1}bd^qc^s)^{-q-1}(dc^{-s}d^{-q-1}ad^qc^s)(dc^{-s}d^{-q-1}bd^qc^s)^q.
\end{equation}
If $q=-1$ (case 3.1.1.1), then  equality (\ref{equality2}) clearly holds, i.~e. for $W(x,y,z)=yz^{-s}xy^{-1}z^s$ the algebraic system $(G,*_{W,c})$ is a quandle for every group $G$ and an element $c\in G$. It means that in Case 3.1.1.1 we have a quandle. If $q\neq 1$ (case 3.1.1.2), i.~e. $q<-1$ or $q>0$, then look to $a,b,c,d$ as to the generators of the free group $F(a,b,c,d)$. The reduced word on the left side of equality (\ref{equality2}) begins with $c^{-s}$, while the word on the right side of equality (\ref{equality2}) begins either with $d^{q+1}$ (if $q<-1$) or with $d^q$ (of $q>0$), hence equality (\ref{equality2}) doesn't hold. It means that in Case 3.1.1.2 we don't have a quandle.

Let us consider the case when $\varepsilon=-1$ (case 3.1.2). From equalities (\ref{left3}), (\ref{right3}) we see that the equality $(a*_{W,c}b)*_{W,c}d=(a*_{W,c}d)*_{W,c}(b*_{W,c}d)$ holds if an only if 
\begin{multline}\label{equality22}
c^{-s}b^{-q}ab^{q-1}c^sb^{-1}d^q=\\
=b^{-1}d^q(dc^{-s}d^{-q+1}b^{-1}d^qc^s)^{-q+1}c^{-s}d^{-q}ad^{q-1}c^sd^{-1}(dc^{-s}d^{-q+1}b^{-1}d^qc^s)^q.
\end{multline}
If $q=1$ (case 3.1.2.1), then equality (\ref{equality22}) can be rewritten in the following form
$$
	c^{-s}b^{-1}ac^sb^{-1}d=
	b^{-1}dc^{-s}d^{-1}ab^{-1}dc^s.
$$
Looking to $a,b,c,d$ as to the generators of the free group it is clear that the last equality doesn't hold for all non-zero integers $s$, hence, if $q=1$, then  there exists a group $G$ and an element $c\in G$ such that $(G,*_{W,c})$ is not a quandle. It means that in Case 3.1.2.1 we don't have a quandle. If $q\neq 1$, i.~e. $q<0$ or $q>1$ then look to $a,b,c,d$ as to the generators of the free group. The reduced word on the left side of equality (\ref{equality22}) begins with $c^{-s}$, while the word on the right side of equality (\ref{equality22}) begins with $b^{-1}$, hence, equality (\ref{equality22}) doesn't hold. So, in this case there exists a group $G$ and an element $c\in G$ such that $(G,*_{W,c})$ is not a quandle. It means that in Case 3.1.2.2 we don't have a quandle.

Let us now consider the case when $|w(y,z)|=2$ and $w(y,z)=z^sy^q$  for non-zero integers $q,s\in\mathbb{Z}$ (case 3.2). In this case the word $W(x,y,z)$ has the following form 
$$W(x,y,z)=y^{1-q}z^{-s}y^{-\varepsilon}x^{\varepsilon}z^sy^q.$$
Let us understand when the operation $*_{W,c}$ satisfies the third quandle axiom (q3) for each group $G$ and an element $c\in G$
\begin{align}
	\notag(a*_{W,c}b)*_{W,c}d&=(W(a,b,c))*_{W,c}d\\
	\notag&=(b^{1-q}c^{-s}b^{-\varepsilon}a^{\varepsilon}c^sb^q)*_{W,c}d\\
	\notag&=W(b^{1-q}c^{-s}b^{-\varepsilon}a^{\varepsilon}c^sb^q,d,c)\\
	\label{left33}&=d^{1-q}c^{-s}d^{-\varepsilon}\left(b^{1-q}c^{-s}b^{-\varepsilon}a^{\varepsilon}c^sb^q\right)^{\varepsilon}c^sd^q,\\
	\notag 	(a*_{W,c}d)*_{W,c}(b*_{W,c}d)&=(W(a,d,c))*_{W,c}(W(b,d,c))\\
	\notag&=(d^{1-q}c^{-s}d^{-\varepsilon}a^{\varepsilon}c^sd^q)*_{W,c}(d^{1-q}c^{-s}d^{-\varepsilon}b^{\varepsilon}c^sd^q)\\
	\notag&=W(d^{1-q}c^{-s}d^{-\varepsilon}a^{\varepsilon}c^sd^q,d^{1-q}c^{-s}d^{-\varepsilon}b^{\varepsilon}c^sd^q,c)\\
\notag	&=\left(d^{1-q}c^{-s}d^{-\varepsilon}b^{\varepsilon}c^sd^q\right)^{1-q}c^{-s}\left(d^{1-q}c^{-s}d^{-\varepsilon}b^{\varepsilon}c^sd^q\right)^{-\varepsilon}\\
	\label{right33}&~~~~\cdot\left(d^{1-q}c^{-s}d^{-\varepsilon}a^{\varepsilon}c^sd^q\right)^{\varepsilon}c^s\left(d^{1-q}c^{-s}d^{-\varepsilon}b^{\varepsilon}c^sd^q\right)^q
	\end{align}
Consider the case when $\varepsilon=1$ (case 3.2.1). From equalities (\ref{left33}), (\ref{right33}) we see that the equality $(a*_{W,c}b)*_{W,c}d=(a*_{W,c}d)*_{W,c}(b*_{W,c}d)$ holds if an only if 
\begin{multline}\label{longerthing}
	d^{1-q}c^{-s}d^{-1}b^{1-q}c^{-s}b^{-1}ac^sb^qc^sd^q=\\
	=\left(d^{1-q}c^{-s}d^{-1}bc^sd^q\right)^{1-q}c^{-s}d^{-q}c^{-s}b^{-1}ac^sd^qc^s\left(d^{1-q}c^{-s}d^{-}bc^sd^q\right)^q.
\end{multline}
If $q=1$ (case 3.2.1.1), then equality (\ref{longerthing}) clearly holds independently on $s$. It means that in Case 3.2.1.1 we don't have a quandle. If $q\neq 1$ (case 3.2.1.2), i.~e. $q<0$ or $q>1$, then look to $a,b,c,d$ as to the generators of the free group. If $q<0$, then the reduced word on the left side of equality (\ref{longerthing}) ends with $d^{q}$, while the word on the right side of equality (\ref{longerthing}) ends with $d^{q-1}$, hence, equality (\ref{longerthing}) doesn't hold. If $q>1$, then the reduced word on the left side of equality (\ref{longerthing}) begins with $d^{1-q}$, while the word on the right side of equality (\ref{longerthing}) begins with $d^{-q}$, hence, equality (\ref{longerthing}) doesn't hold. It means that in Case 3.2.1.2 we don't have a quandle.

Let us consider the case when $\varepsilon=-1$ (case 3.2.2). From equalities (\ref{left33}), (\ref{right33}) we see that the equality $(a*_{W,c}b)*_{W,c}d=(a*_{W,c}d)*_{W,c}(b*_{W,c}d)$ holds if an only if 
\begin{multline}\label{anotherlong}
d^{1-q}c^{-s}db^{-q}c^{-s}ab^{-1}c^sb^{q-1}c^sd^q=\\
		=\left(d^{1-q}c^{-s}db^{-1}c^sd^q\right)^{1-q}c^{-s}d^{1-q}c^{-s}db^{-1} ad^{-1}c^sd^{q-1}c^s\left(d^{1-q}c^{-s}db^{-1}c^sd^q\right)^q.
\end{multline}
Look to $a,b,c,d$ as to the generators of the free group. If $q<0$, then the reduced word on the left side of equality (\ref{anotherlong}) ends with $d^{q}$, while the word on the right side of equality (\ref{anotherlong}) ends with $d^{q-1}$, hence, equality (\ref{anotherlong}) doesn't hold. If $q=1$, then using direct calculations it is clear that equality (\ref{anotherlong}) doesn't hold.  If $q>1$, then the reduced word on the left side of equality (\ref{anotherlong}) begins with $d^{1-q}$, while the word on the right side of equality (\ref{anotherlong}) begins with $d^{-q}$, hence, equality (\ref{anotherlong}) doesn't hold. It means that in Case 3.2.2 we don't have a quandle.

\textbf{Case 4:} $|w(y,z)|=3$. In this case $w(y,z)$ is either $w(y,z)=y^{q_1}z^{s}y^{q_2}$ for non-zero integers $q_1,q_2,s$ or $w(y,z)=z^{s_1}y^{q}z^{s_2}$ for non-zero integers $q,s_1,s_2$. 

Consider the case when $w(y,z)=y^{q_1}z^{s}y^{q_2}$ for non-zero integers $q_1,q_2,s$ (case~4.1). In this case $W(x,y,z)=y\left(y^{q_1}z^{s}y^{q_2}\right)^{-1}y^{-\varepsilon}x^{\varepsilon}\left(y^{q_1}z^{s}y^{q_2}\right)=y^{1-q_2}z^{-s}y^{-\varepsilon-q_1}x^{\varepsilon}y^{q_1}z^sy^{q_2}$.  
Let us understand when the operation $*_{W,c}$ satisfies the third quandle axiom (q3) for each group $G$ and an element $c\in G$
  \begin{align}
	\notag(a&*_{W,c}b)*_{W,c}d=\\
	\notag&=(W(a,b,c))*_{W,c}d\\
	\notag&=(b^{1-q_2}c^{-s}b^{-\varepsilon-q_1}a^{\varepsilon}b^{q_1}c^sb^{q_2})*_{W,c}d\\
	\notag&=W\left(b^{1-q_2}c^{-s}b^{-\varepsilon-q_1}a^{\varepsilon}b^{q_1}c^sb^{q_2},d,c\right)\\
	\label{left7}&=d^{1-q_2}c^{-s}d^{-\varepsilon-q_1}\left(b^{1-q_2}c^{-s}b^{-\varepsilon-q_1}a^{\varepsilon}b^{q_1}c^sb^{q_2}\right)^{\varepsilon}d^{q_1}c^sd^{q_2}\\
	\notag 	(a&*_{W,c}d)*_{W,c}(b*_{W,c}d)=\\
	\notag&=(W(a,d,c))*_{W,c}(W(b,d,c))\\
	\notag&=\left(d^{1-q_2}c^{-s}d^{-\varepsilon-q_1}a^{\varepsilon}d^{q_1}c^sd^{q_2}\right)*_{W,c}\left(d^{1-q_2}c^{-s}d^{-\varepsilon-q_1}b^{\varepsilon}d^{q_1}c^sd^{q_2}\right)\\
	\notag&=W(d^{1-q_2}c^{-s}d^{-\varepsilon-q_1}a^{\varepsilon}d^{q_1}c^sd^{q_2}, d^{1-q_2}c^{-s}d^{-\varepsilon-q_1}b^{\varepsilon}d^{q_1}c^sd^{q_2},c)\\
	\notag&=\left(d^{1-q_2}c^{-s}d^{-\varepsilon-q_1}b^{\varepsilon}d^{q_1}c^sd^{q_2}\right)^{1-q_2}c^{-s}\left(d^{1-q_2}c^{-s}d^{-\varepsilon-q_1}b^{\varepsilon}d^{q_1}c^sd^{q_2}\right)^{-\varepsilon-q_1}\\
	\notag&~~~~\cdot\left(d^{1-q_2}c^{-s}d^{-\varepsilon-q_1}a^{\varepsilon}d^{q_1}c^sd^{q_2}\right)^{\varepsilon}\left(d^{1-q_2}c^{-s}d^{-\varepsilon-q_1}b^{\varepsilon}d^{q_1}c^sd^{q_2}\right)^{q_1}\\
	\label{right7}&~~~~\cdot c^s\left(d^{1-q_2}c^{-s}d^{-\varepsilon-q_1}b^{\varepsilon}d^{q_1}c^sd^{q_2}\right)^{q_2}
\end{align}
Consider the case when $\varepsilon=1$ (case 4.1.1). From equalities (\ref{left7}), (\ref{right7}) we see that the equality $(a*_{W,c}b)*_{W,c}d=(a*_{W,c}d)*_{W,c}(b*_{W,c}d)$ holds if an only if 
\begin{align}
	\notag d^{1-q_2}c^{-s}&d^{-1-q_1}b^{1-q_2}c^{-s}b^{-1-q_1}ab^{q_1}c^sb^{q_2}d^{q_1}c^sd^{q_2}=\\
	\notag&=\left(d^{1-q_2}c^{-s}d^{-1-q_1}bd^{q_1}c^sd^{q_2}\right)^{1-q_2}c^{-s}\left(d^{1-q_2}c^{-s}d^{-1-q_1}bd^{q_1}c^sd^{q_2}\right)^{-1-q_1}\\
	\notag&~~~~\cdot d^{1-q_2}c^{-s}d^{-1-q_1}ad^{q_1}c^sd^{q_2}\left(d^{1-q_2}c^{-s}d^{-1-q_1}bd^{q_1}c^sd^{q_2}\right)^{q_1}\\
	\label{q2=1}&~~~~\cdot c^s\left(d^{1-q_2}c^{-s}d^{-1-q_1}bd^{q_1}c^sd^{q_2}\right)^{q_2}
\end{align}
If $q_2=1$, then equality (\ref{q2=1}) can be rewritten in the following form
\begin{multline}\label{q2=11}
	d^{-1-q_1}c^{-s}b^{-1-q_1}ab^{q_1}c^s=\\
	=\left(c^{-s}d^{-1-q_1}bd^{q_1}c^sd\right)^{-1-q_1}c^{-s}d^{-1-q_1}ad^{q_1}c^sd\left(c^{-s}d^{-1-q_1}bd^{q_1}c^sd\right)^{q_1}d^{-1-q_1}
\end{multline}
If $q_1=-1$ (case 4.1.1.1), then equality (\ref{q2=11}) holds independently on $s$. It means that in Case 4.1.1.1 we have a quandle. If $q_1\neq 1$ (case 4.1.1.2), i.~e. $q_1<-1$ or $q_1>0$, then look to $a,b,c,d$ as to the generators of the free group. If $q_1<-1$, then the reduced word on the left side of equality (\ref{q2=11}) end with $c^{s}$, while the word on the right side of equality (\ref{q2=11}) end with $d^{-1-q_1}$, hence, equality (\ref{q2=11}) doesn't hold. If $q_1>0$, then the reduced word on the left side of equality (\ref{q2=11}) ends with $c^s$, while the word on the right side of equality (\ref{q2=11}) ends with $d^{-q_1}$, hence, equality (\ref{q2=11}) doesn't hold. It means that in Case 4.1.1.2 we don't have a quandle.

If $q_2\neq1$ (case 4.1.1.3), i.~e. $q_2<0$ or $q_2>1$, then look to $a,b,c,d$ as to the generators of the free group. If $q_2<0$, then the word on the left side of equality (\ref{q2=1}) ends with $d^{q_2}$, while the word on the right side of equality (\ref{q2=1}) ends with $d^{q_2-1}$, hence, equality (\ref{q2=1}) doesn't hold. If $q_2>1$, then the word on the left side of equality (\ref{q2=1}) begins with $d^{1-q_2}$, while the word on the right side of equality (\ref{q2=1}) begins with $d^{-q_2}$, hence, equality (\ref{q2=1}) doesn't hold. It means that in Case 4.1.1.3 we don't have a quandle.

Codsider the case when $\varepsilon=-1$ (case 4.1.2). From equalities (\ref{left7}), (\ref{right7}) we see that the equality $(a*_{W,c}b)*_{W,c}d=(a*_{W,c}d)*_{W,c}(b*_{W,c}d)$ holds if an only if 
\begin{align}
\notag 	d^{1-q_2}c^{-s}&d^{1-q_1}b^{-q_2}c^{-s}b^{-q_1}ab^{q_1-1}c^sb^{q_2-1}d^{q_1}c^sd^{q_2}=\\
	\notag&=\left(d^{1-q_2}c^{-s}d^{1-q_1}b^{-1}d^{q_1}c^sd^{q_2}\right)^{1-q_2}c^{-s}\left(d^{1-q_2}c^{-s}d^{1-q_1}b^{-1}d^{q_1}c^sd^{q_2}\right)^{1-q_1}\\
	\notag&~~~~\cdot d^{-q_2}c^{-s}d^{-q_1}ad^{q_1-1}c^sd^{q_2-1}\left(d^{1-q_2}c^{-s}d^{1-q_1}b^{-1}d^{q_1}c^sd^{q_2}\right)^{q_1}\\
\label{e-1last}	&~~~~\cdot c^s\left(d^{1-q_2}c^{-s}d^{1-q_1}b^{-1}d^{q_1}c^sd^{q_2}\right)^{q_2}
\end{align}
If $q_2=1$ (case 4.1.2.1), then equality (\ref{e-1last}) can be rewritten in the following form
\begin{align}
	\notag 	&d^{1-q_1}b^{-1}c^{-s}b^{-q_1}ab^{q_1-1}c^s=\\
\label{e-1last2}	&=\left(c^{-s}d^{1-q_1}b^{-1}d^{q_1}c^sd\right)^{1-q_1}d^{-1}c^{-s}d^{-q_1}ad^{q_1-1}c^s\left(c^{-s}d^{1-q_1}b^{-1}d^{q_1}c^sd\right)^{q_1}d^{1-q_1}b^{-1}
\end{align}
If we look to $a,b,c,d$ as to the generators of the free group, then it is clear that the word on the left side of equality (\ref{e-1last2}) ends with $c^{s}$, while the word on the right side of equality (\ref{e-1last2}) ends with $b^{-1}$, hence, equality (\ref{e-1last2}) doesn't hold. It means that in Case 4.1.2.1 we don't have a quandle.

If $q_2\neq1$ (case 4.1.2.2), i.~e. $q_2<0$ or $q_2>1$, then look to $a,b,c,d$ as to the generators of the free group. If $q_2<0$, then the word on the left side of equality (\ref{e-1last}) ends with $d^{q_2}$, while the word on the right side of equality (\ref{e-1last}) ends with $d^{q_2-1}$, hence, equality (\ref{e-1last}) doesn't hold. If $q_2>1$, then the word on the left side of equality (\ref{e-1last}) begins with $d^{1-q_2}$, while the word on the right side of equality (\ref{e-1last}) begins with $d^{-q_2}$, hence, equality (\ref{e-1last}) doesn't hold. It means that in Case 4.1.2.2 we don't have a quandle.

Let us consider the case when $w(y,z)=z^{s_1}y^{q}z^{s_2}$ for non-zero integers $q,s_1,s_2$ (case~4.2). In this case 
\begin{equation}\label{generalbigguy}
W(x,y,z)=y\left(z^{s_1}y^{q}z^{s_2}\right)^{-1}y^{-\varepsilon}x^{\varepsilon}\left(z^{s_1}y^{q}z^{s_2}\right)=yz^{-s_2}y^{-q}z^{-s_1}y^{-\varepsilon}x^{\varepsilon}z^{s_1}y^{q}z^{s_2}.
\end{equation}
Let us understand when the operation $*_{W,c}$ satisfies the third quandle axiom (q3) for each group $G$ and an element $c\in G$
 \begin{align}
	\notag(a*_{W,c}b)*_{W,c}d&=(W(a,b,c))*_{W,c}d\\
	\label{leftpartalmostfinal}&=W(W(a,b,c),d,c)\\
	\notag(a*_{W,c}d)*_{W,c}(b*_{W,c}d)&=(W(a,d,c))*_{W,c}(W(b,d,c))\\
	\label{righttpartalmostfinal}&=W(W(a,d,c),W(b,d,c),c)
\end{align}
Let us look to the elements $a,b,c,d$ as to the generators of the free group. For a word $A\in F(a,b,c,d)$ denote by $|A|_c$ the number of syllables which are the powers of the element $c$ in the reduced form of $A$. Let us evaluate the values $|(a*_{W,c}b)*_{W,c}d|_c$ and $|(a*_{W,c}d)*_{W,c}(b*_{W,c}d)|_c$ given by formulas (\ref{leftpartalmostfinal}), (\ref{righttpartalmostfinal}).

First, let us find the upper bound for the value $|(a*_{W,c}b)*_{W,c}d|_c$. Denote by $A_1=W(a,b,c)$. From equality (\ref{leftpartalmostfinal}) and equality (\ref{generalbigguy}) it follows that 
 \begin{align}
	\label{leftpartevaluation}(a*_{W,c}b)*_{W,c}d&=W(A,d,c)=dc^{-s_2}d^{-q}c^{-s_1}d^{-\varepsilon}A_1^{\varepsilon}c^{s_1}d^{q}c^{s_2}
\end{align}
In order to calculate the upper bound for $|(a*_{W,c}b)*_{W,c}d|_c$ we have to calculate the number of syllables which are the powers of the element $c$ in the right part of equality (\ref{leftpartevaluation}) without thinking if some powers of $c$ can be reduced with each other (since we try to calculate the upper bound). We have the following upper bound 
 \begin{align}
	\notag |(a*_{W,c}b)*_{W,c}d|_c\leq 1+1+|A_1|+1+1=|A_1|+4.
\end{align}

Let us now calculate the lower bound for the value $|(a*_{W,c}d)*_{W,c}(b*_{W,c}d)|_c$. Denote by $A_2=W(a,d,c)$, $B=W(b,d,c)$. From equality (\ref{righttpartalmostfinal}) and equality (\ref{generalbigguy}) it follows that 
 \begin{align}
	\label{rightpartevaluation}(a*_{W,c}d)*_{W,c}(b*_{W,c}d)&=W(A_2,B,c)=Bc^{-s_2}B^{-q}c^{-s_1}B^{-\varepsilon}A_2^{\varepsilon}c^{s_1}B^{q}c^{s_2}.
\end{align}
In order to calculate the lower bound for $|(a*_{W,c}d)*_{W,c}(b*_{W,c}d)|_c$ we have to calculate the number of syllables which are powers of $c$ in the right part of equality (\ref{rightpartevaluation}) which are guaranteed not to be reduced or merged with other syllables which are powers of $c$. 

From the equality $B=W(b,d,c)$ and equality (\ref{generalbigguy}) it follows that if $|q|\neq1$ (case 4.2.1), then in the word $B^{\pm q}$ there are at least $4(q-1)$ syllables which are powers of $c$  which are guaranteed not to be reduced or merged with other syllables which are powers of $c$ in the right part of equality (\ref{rightpartevaluation}). This is due to the fact that between $2$ consequent syllables which are powers of $b$ in the word $B^{\pm q}=W(b,d,c)^{\pm q}$ there are exactly $4$ syllables which are powers of $c$  which are guaranteed not to be reduced or merged with other syllables which are powers of $c$, and in the word $B^{\pm q}=W(b,d,c)^{\pm q}$ there are exactly $q$ syllables which are powers of $b$, hence, there are at least $4(q-1)$ syllables which are powers of $c$. Due to this remark in order to calculate the lower bound for $|(a*_{W,c}d)*_{W,c}(b*_{W,c}d)|_c$ let in the right part of equality (\ref{rightpartevaluation}) calculate only syllables which are powers of $c$ in the word $B^q$, $B^{-q}$ and $A_2^{\varepsilon}$ (these syllables which are powers of $c$  are guaranteed not to be reduced or merged with other syllables which are powers of $c$). We have the following lower bound
 \begin{align}
	\notag |(a*_{W,c}d)*_{W,c}(b*_{W,c}d)|_c\geq 4(q-1)+|A_2|_c+4(q-1)=|A_2|_c+8(q-1).
\end{align}
Hence we have proved that the following inequalities hold 
 \begin{align}
	\notag |(a*_{W,c}b)*_{W,c}d|_c&\leq |A_1|_c+4,\\
		\notag |(a*_{W,c}d)*_{W,c}(b*_{W,c}d)|_c&\geq |A_2|_c+8(q-1).
\end{align}
Noting that $|A_1|_c=|W(a,b,c)|_c=|W(a,d,c)|_c=|A_2|_c$ we conclude that 
$$|(a*_{W,c}b)*_{W,c}d|_c\neq |(a*_{W,c}d)*_{W,c}(b*_{W,c}d)|_c$$
and hence $(a*_{W,c}b)*_{W,c}d\neq (a*_{W,c}d)*_{W,c}(b*_{W,c}d)$, i.~e. in this case there exists a group $G$ and an element $c\in G$ such that $(G,*_{W,c})$ is not a quandle. It means that in Case 4.2.1 we don't have a quandle.

Let us now consider the case when $|q|=1$ (case 4.2.2), i.~e. the case when either $q=1$, or $q=-1$. In this case we can directly calculate all syllabels which are powers of $c$ in both  $(a*_{W,c}b)*_{W,c}d$ and $(a*_{W,c}d)*_{W,c}(b*_{W,c}d)$, and make sure that $|(a*_{W,c}b)*_{W,c}d|_c\neq |(a*_{W,c}d)*_{W,c}(b*_{W,c}d)|_c$. Hence
$$(a*_{W,c}b)*_{W,c}d\neq (a*_{W,c}d)*_{W,c}(b*_{W,c}d),$$ 
i.~e. in this case there exists a group $G$ and an element $c\in G$ such that $(G,*_{W,c})$ is not a quandle. It means that in Case 4.2.2 we don't have a quandle.

\textbf{Case 5:} $|w(y,z)|\geq4$. In this case the word $w(y,z)$ has at least two syllables which are powers of $z$. Using the same technique of evaluating the values  
\begin{align*}|(a*_{W,c}b)*_{W,c}d|_c&& \text{and}&&|(a*_{W,c}d)*_{W,c}(b*_{W,c}d)|_c
	\end{align*}
as we did in case 4.2.1, repeating the proof of case 4.2.1 almost word by word we conclude that 
\begin{align*}|(a*_{W,c}b)*_{W,c}d|_c<|(a*_{W,c}d)*_{W,c}(b*_{W,c}d)|_c,
\end{align*}
there exists a group $G$ and an element $c\in G$ such that $(G,*_{W,c})$ is not a quandle. It means that in Case 5 we don't have a quandle.\hfill$\square$
\section{Verbal quandles with $n$ parameters}
In this short section, we note that there is no sense to classify verbal quands with more than one parameter. More precisely, using the technique used for proving Theorem~\ref{th22233}, we can show that the following statement about verbal quands with $n$ parameters holds.
\begin{ttt}
Let $W\in F(x,y,z_1,z_2,\dots,z_n)$ be such that $Q=(G,*_{W,c_1,\dots,c_n})$ is a quandle for every group $G$ and elements $c_1,c_2,\dots,c_n\in G$. Then $W$ has one of the following forms
\linespread{0}
\begin{enumerate}
	\item $W=yx^{-1}y$,
	\item $W=y^{-n}xy^n$ for $n\in\mathbb{N}$,
	\item $W=yu^{-1}xy^{-1}u$ for $u\in F(z_1,z_2,\dots,z_n)$,
	\item $W=u^{-1}xy^{-1}uy$ for $u\in F(z_1,z_2,\dots,z_n)$,
	\item $W=yu^{-1}y^{-1}xu$ for $u\in F(z_1,z_2,\dots,z_n)$,
	\item $W=u^{-1}y^{-1}xuy$ for $u\in F(z_1,z_2,\dots,z_n)$.
\end{enumerate} 
\end{ttt}

~\\

\noindent Elizaveta Markhinina\\
Novosibirsk State University, \underline{e.markhinina@g.nsu.ru}

~\\
Timur Nasybullov\\
Novosibirsk State University, Sobolev Institute of Mathematics, \underline{timur.nasybullov@mail.ru}
\end{document}